\documentclass[11pt]{amsart}

\usepackage{amsthm}
\usepackage{amssymb}
\usepackage{amsmath}
\usepackage{amsfonts} 

\newcommand{\excise}[1]{}

\newtheorem{thm}{Theorem}
\newtheorem{lemma}[thm]{Lemma}

\newtheorem{prop}[thm]{Proposition}

\theoremstyle{definition}

\newtheorem{remark}[thm]{Remark}

\newtheorem{defn}[thm]{Definition}
\newtheorem{conv}[thm]{Convention}

        {\begin{list}
                {\noindent\makebox[0mm][r]{\arabic{enumi}.}}
                {\leftmargin=5.5ex \usecounter{enumi}}
        }
        {\end{list}}

        {\begin{list}
                {\noindent\makebox[0mm][r]{(\roman{enumi})}}
                {\leftmargin=5.5ex \usecounter{enumi}}
        }
        {\end{list}}

\def\bem#1{\textbf{#1}}

\def\<{\langle}
\def\>{\rangle}
\def\CC{{\mathbb C}}
\def\EE{{\mathcal E}}
\def\NN{{\mathbb N}}
\def\RR{{\mathbb R}}
\def\ZZ{{\mathbb Z}}
\def\mm{{\mathfrak m}}

\def\th{{\rm th}}

\def\vol{{\rm vol}}
\def\conv{{\rm conv}}
\def\rank{{\rm rank}}
\def\cocoa
  {{\hbox{\rm C\kern-.13em o\kern-.07em C\kern-.13em o\kern-.15em A}}}

\def\minus{\smallsetminus}
\def\implies{\Rightarrow}

\def\ol#1{{\overline {#1}}}

\begin{document}

\mbox{}
\vspace{-3.2ex}
\title{Combinatorics of rank jumps in simplicial hypergeometric systems}
\author{Laura Felicia Matusevich}
\address[Laura Felicia Matusevich]{Mathematical Sciences Research
	Institute\\Berkeley, CA}
\curraddr{Department of Mathematics\\Harvard University\\Cambridge, MA}
\email{laura@math.harvard.edu}
\author{Ezra Miller}
\address[Ezra Miller]{Mathematical Sciences Research Institute\\Berkeley, CA}
\curraddr{School of Mathematics, University of Minnesota,
	Minneapolis MN}
\email{ezra@math.umn.edu}
\date{27 January 2004}

\begin{abstract}
Let $A$ be an integer $d \times n$ matrix, and assume that the convex
hull $\conv(A)$ of its columns is a simplex of dimension~\mbox{$d-1$}.
It is known that the semigroup ring $\CC[\NN A]$ is Cohen--Macaulay if
and only if the rank of the GKZ hypergeometric system $H_A(\beta)$
equals the normalized volume of $\conv(A)$ for all complex parameters
$\beta \in \CC^d$ \cite{mutsumi2}.
Our refinement here shows that $H_A(\beta)$ has rank strictly larger
than the volume of $\conv(A)$ if and only if $\beta$ lies in the
Zariski closure (in~$\CC^d$) of all $\ZZ^d$-graded degrees where the
local cohomology \mbox{$\bigoplus_{i < d} H^i_\mm(\CC[\NN A])$ is
nonzero}.
\end{abstract}
\maketitle

\section{Introduction}

Gelfand, Graev, Kapranov, and Zelevinsky \cite{GGZ,GKZ} defined
certain linear systems of partial differential equations, now known as
{\em $A$-hypergeometric} or {\em GKZ hypergeometric}\/ systems
$H_A(\beta)$, whose solutions generalize the classical hypergeometric
series.  These holonomic systems are constructed from discrete input
consisting of an integer $d \times n$ matrix~$A$, along with
continuous input consisting of a complex vector~\mbox{$\beta \in
\CC^d$}.  The matrix~$A$ defines a semigroup ring $\CC[\NN A]$, and
the same authors have shown that the dimension~$\rank(H_A(\beta))$ of
the space of analytic solutions of~$H_A(\beta)$ is independent
of~$\beta$ when $\CC[\NN A]$ is Cohen--Macaulay \cite{GKZ}.

Meanwhile, Adolphson showed that even when $\CC[\NN A]$ is not
Cohen--Macaulay, the rank of $H_A(\beta)$ is independent of~$\beta$,
as long as $\beta$ is generic in a certain precise sense
\cite{adolphson}.  After Sturmfels and Takayama showed that the rank
can actually go up for non-generic parameters~$\beta$ \cite{StuTaka},
Cattani, D'Andrea and Dickenstein \cite{monomialcurves} showed that if
$\conv(\NN A)$ is a segment, then in fact the rank does jump whenever
$\CC[\NN A]$ fails to be Cohen--Macaulay.  This result was generalized
by Saito \cite{mutsumi2}, who, using different methods, proved that
there existed rank-jumping parameters for any non Cohen--Macaulay
simplicial semigroup $\CC[\NN A]$.

In this note, we use the combinatorics of $\ZZ^d$-graded local
cohomology to \mbox{characterize} the set of parameters~$\beta$ for
which the rank goes up, in the simplicial case.  Our premise, reviewed
in Section~\ref{localcoh}, is the standard fact that a semigroup
ring~$\CC[\NN A]$ fails to be Cohen--Macaulay if and only if a local
cohomology module $H^i_\mm(\CC[\NN A])$ is nonzero for some
cohomological index~$i$ strictly less than than the dimension~$d$
of~$\CC[\NN A]$.  After gathering some facts about $A$-hypergeometric
systems in Section~\ref{hyper}, we show in
Theorem \ref{t:simplicial} that the set of {\em rank jumping
parameters}\/ is the Zariski closure (in~$\CC^d$) of the set of all
$\ZZ^d$-graded degrees where the local cohomology \mbox{$\bigoplus_{i
< d} H^i_\mm(\CC[\NN A])$ is nonzero}.

The original role of this result was as evidence for our conjecture
that it generalizes to arbitrary integer matrices~$A$, regardless
of whether or not the semigroup $\NN A$ is simplicial.  Although we
do not know whether the methods of this note extend to the general
case, our conjecture has since been proved using a different
approach~\cite{ranks}.

\section{Computing local cohomology for semigroup rings}\label{localcoh}

Throughout this note, let $A$ be a $d\times n$ integer matrix whose
first row has all entries equal to~$1$, and whose columns
$a_1,\ldots,a_n$ generate $\ZZ^d$ as a group.  Unless otherwise
explicitly stated, we do not assume that the polytope $\conv(A)$
obtained by taking the convex hull (in~$\RR^d$) of the column
vectors $a_1,\ldots,a_n$ is a simplex; in particular, we need no
simplicial assumptions from here through Definition~\ref{d:excep}.
The semigroup
\begin{eqnarray*}
  \NN A &=& \Big\{\sum_{i=1}^n k_i a_i \mid k_1,\dots ,k_n \in \NN\Big\}
\end{eqnarray*}
has \bem{semigroup ring}\/ $R=\CC[\NN A] \cong
\CC[\partial_1,\ldots,\partial_n]/I_A$, where
\begin{eqnarray*}
  I_A &=& \<\partial^u-\partial^v \mid A\cdot u = A\cdot v\>
\end{eqnarray*}
is the \bem{toric ideal} of~$A$.  The ring~$R$ is naturally graded
by~$\ZZ^d$, with the $i^\th$ indeterminate having degree
$\deg(\partial_i) = a_i$ equal to the $i^\th$ column of~$A$.
By~a \bem{face} of~$\NN A$ we mean a set of lattice points minimizing
some linear functional on~$\NN A$.  The terms \bem{ray}
and~\bem{facet} refer to faces of dimension~$1$ and~\mbox{$d-1$},
respectively, where the dimension of a face equals the rank of the
subgroup $\ZZ \tau \subseteq \ZZ^d$ it generates.  It is convenient to
identify a face~$\tau$ of~$\NN A$ with the subset of~$\{1,\ldots,n\}$
indexing the \mbox{vectors~$a_i$ lying in~$\tau$}.

We now recall some facts from \cite[Chapter~6]{BH} about the local
cohomology modules~$H^i_\mm(R)$, where $\mm = \<\partial_1,\dots
\partial_n\>$ is the graded maximal ideal of~$R$.  Since $R$ is a
semigroup ring, the local cohomology of~$R$ is the cohomology of the
complex
\begin{equation} \label{eqn:complex}
  0 \rightarrow R \to
  \bigoplus_{\text{rays } \tau} R_\tau \to
  \bigoplus_{\text{2-dim faces }\tau} R_\tau \to
  \cdots \to
  \bigoplus_{\text{facets }\tau} R_\tau
  \to R_\mm \to 0,
\end{equation}
where $R_\tau$ is the localization of~$R$ by inverting the
indeterminates $\partial_i$ for $i \in \tau$.  The differential is
derived from the algebraic cochain complex of the polytope~$\conv(A)$,
once orientations on the faces of~$\conv(A)$ have been chosen.

The above local cohomology can be computed multidegree by multidegree.
Indeed, the localization $R_\tau$ is nonzero in graded degree $\beta
\in \ZZ^d$ if and only if $\beta$ lies in the subsemigroup $\NN A +
\ZZ \tau$ of~$\ZZ^d$; in other words, $R_\tau = \CC[\NN A + \ZZ\tau]$.
Therefore, the faces of~$\NN A$ contributing a nonzero vector space
(of dimension~$1$) to the degree~$\beta$ piece of the
complex~\eqref{eqn:complex} is
\begin{eqnarray*}
  \nabla(\beta) &=& \{\text{faces }\tau \text{ of }\NN A \mid \beta
  \in \NN A + \ZZ \tau\}.
\end{eqnarray*}
This set of faces is closed under going up, meaning that if $\tau
\subset \sigma$ and $\tau \in \nabla(\beta)$, then also $\sigma \in
\nabla(\beta)$.  When we write cohomology groups $H^j(\nabla)$ for such a
\bem{polyhedral cocomplex}, what we mean formally is that
\begin{eqnarray*}
  H^j(\nabla) &=& H^j(\conv(A), \conv(A)\!\minus\!\nabla; \CC)
\end{eqnarray*}
is the cohomology with complex coefficients of the relative cochain
complex of the complementary polyhedral subcomplex of~$\conv(A)$.

Here is a standard result in combinatorial commutative algebra.

\begin{thm}
The local cohomology $H^j_\mm(R)_\beta$ of the semigroup ring~$R$ in
$\ZZ^d$-graded degree~$\beta$ is isomorphic to~$H^j(\nabla(\beta))$.
In particular, $R$ is Cohen--Macaulay if and only if
$H^j(\nabla(\beta)) = 0$ for all degrees $\beta \in \ZZ^d$ and
cohomological degerees $j=0,\ldots,d-1$.%
\end{thm}
\begin{proof}
Use the complex in~\eqref{eqn:complex} to compute local cohomology.
\end{proof}

\begin{defn}{\textbf{\cite{injAlg}}} \label{d:sector}
The \bem{sector partition} is the partition of $\ZZ^d$ into
equivalence classes for which $\beta \equiv \beta'$ if and only if
$\nabla(\beta) = \nabla(\beta')$.  For a cocomplex~$\nabla$,
the (possibly empty) of degrees $\beta \in \ZZ^d$ satisfying
$\nabla(\beta) = \nabla$ is a \bem{sector}.
\end{defn}

Since the local cohomology of~$R$ in degree $\beta$ only depends
on~$\nabla(\beta)$, it is constant on every sector.

\begin{defn} \label{d:exdegree}
A degree $\beta \in \ZZ^d$ such that $H^j(\nabla(\beta))\neq 0$ for
some $0 \leq j \leq d-1$ is called an \bem{exceptional degree} of~$A$.
The Zariski closure of the set $E(A)$ of exceptional degrees of~$A$ is
called the \bem{slab arrangement} $\ol E(A)$ of~$A$.  Given an
irreducible component of the slab arrangement, the set of exceptional
degrees lying inside that component and in no other components is
called a~\bem{slab}.%
\end{defn}

\begin{prop} \label{p:slab}
The slab arrangement is a union of affine translates of linear
subspaces $\CC \tau$ generated by faces~$\tau$ of\/~$\NN A$, thought
of as subsets of\/~$\CC^d$.%
\end{prop}
\begin{proof}
The Matlis dual of each local cohomology module is finitely generated,
and therefore has a filtration whose successive quotients are
$\ZZ^d$-graded shifts of quotients of~$R$ by prime monomial ideals.
Each successive quotient is therefore a $\ZZ^d$-graded shift of a
semigroup ring~$\CC[\tau]$ for some face $\tau$ of~$\NN A$.  The
Matlis dual of local cohomology is thus supported on a set of degrees
satisfying the conclusion of the proposition.  The exceptional degrees
are the negatives of the support degrees of the Matlis dual.%
\end{proof}

\section{Local cohomology and \(A\)-hypergeometric systems}\label{hyper}

Denote by $D_n$ the \bem{Weyl algebra}, by which we mean the ring of
linear partial differential operators with polynomial coefficients in
$n$ variables.  That is, $D_n$ is the free associative algebra $\CC
\langle x_1,\dots,x_n, \partial_1,\dots,\partial_n\rangle$ modulo the
relations $x_ix_j-x_jx_i$, $\partial_i\partial_j-\partial_j\partial_i$
and $\partial_j x_i -x_i\partial_j-\delta_{ij}$, where $\delta_{ij}$
is the Kronecker delta.

\begin{defn}
Given $A=(a_{ij})$ as before and $\beta \in \CC^d$, the
\bem{$A$-hypergeometric system} with \bem{parameter~$\beta$} is the
left ideal in the Weyl algebra $D_n$ generated by
\begin{eqnarray*}
  I_A &\text{and}& \sum_{j=1}^n a_{ij}x_j \partial_j -\beta_i \quad
  \text{for } i=1,\dots,d.
\end{eqnarray*}
The \bem{$A$-hypergeometric module} with \bem{parameter~$\beta$} is
$M_A(\beta) = D_n/\!H_A(\beta)$.
\end{defn}

The following result relates our way of computing local cohomology of
semigroup rings to $A$-hypergeometric systems.

\begin{thm}
\label{thm:refinement}
Stratify $\ZZ^d$ so that $\beta$ lies in the same stratum as~$\beta'$
iff the $D$-modules $M_A(\beta)$ and~$M_A(\beta')$ are isomorphic.
This stratification refines the sector partition, meaning that
$M_A(\beta) \cong M_A(\beta')$ implies $\nabla(\beta) =
\nabla(\beta')$.
\end{thm}

To prove this result, we recall Saito's combinatorial results on
isomorphisms of hypergeometric $D$-modules.

\begin{defn}
Let $\beta \in \CC^d$ and $\tau$ a face of the cone $\NN A$.  Let
\begin{eqnarray*}
  E_{\tau}(\beta) &=& \{\lambda \in \CC\tau \mid \beta \in \lambda +
  \NN A + \ZZ \tau\} / \ZZ\tau
\end{eqnarray*}
be the set of vectors $\lambda \in \CC\tau$, up to translation
by~$\ZZ\tau$, for which $\beta - \lambda$ lies in the localization of
$\NN A$ along~$\tau$.
\end{defn}

\begin{thm}\cite{mutsumi2} \label{thm:dmoduleisom}
The $D$-modules $M_A(\beta)$ and $M_A(\beta')$ are isomorphic for two
parameters $\beta$ and~$\beta'$ in~$\CC^d$ if and only if
$E_{\tau}(\beta)=E_{\tau}(\beta')$ for all faces~$\tau$ of~$\NN A$.
\end{thm}

\begin{proof}[Proof of Theorem \ref{thm:refinement}]
Given a vector $\beta \in \ZZ^d$, we have
\begin{eqnarray} \label{nabla}
  \nabla(\beta) &=& \{\text{faces }\tau\text{ of } \NN A \mid 0 \in
  E_\tau(\beta)\}
\end{eqnarray}
by definition.  Therefore, for any pair of parameters $\beta, \beta'
\in \ZZ^d$ such that $M_A(\beta)$ is isomorphic to~$M_A(\beta')$, we
conclude that $\nabla(\beta) = \nabla(\beta')$ by Theorem
\ref{thm:dmoduleisom}.
\end{proof}

\begin{remark}
In general, the refinement in Theorem~\ref{thm:refinement} is proper.
\end{remark}

\section{Rank jumps in the simplex case}\label{rankjump}

\begin{defn} \label{d:excep}
A parameter vector $\beta \in \CC^d$ is an \bem{rank-jumping
parameter} of~$A$ if $\rank(H_A(\beta)) > \vol(A)$, where $\vol(A)$ is
the normalized volume of the polytope $\conv(A)$.  The set of
rank-jumping parameters of~$A$ is called the \bem{exceptional set}
of~$A$, and denoted~$\EE(A)$.
\end{defn}

\begin{thm} \label{t:simplicial}
Fix a $d \times n$ integer matrix~$A$.  If\/ $\conv(A)$ is a
\mbox{$(d-1)$}-simplex, then the exceptional set $\EE(A)$ equals the
Zariski closure~$\ol E(A)$ of the set of exceptional degrees.
\end{thm}

\begin{remark}
Computational evidence (using the computer algebra systems Macaulay~2
\cite{M2}, Singular \cite{Singular}, and CoCoA \cite{cocoa}) as well
as heuristic arguments led us to conjecture the statement of
Theorem~\ref{t:simplicial}.  In fact, the evidence suggested that
Theorem~\ref{t:simplicial} generalizes to the case where $A$ is an
arbitrary integer matrix.
This has since been shown in a subsequent paper \cite{ranks} via
general geometric and homological methods.
\end{remark}

Before getting to the proof, we need four preliminary results.  The
first two do not invoke the hypothesis that $\conv(A)$ is a simplex.

\begin{lemma} \label{l:m}
Suppose that $\rho$ is a face of\/~$\NN A$, and $\alpha \in \rho$ is a
vector not lying on any proper face of~$\rho$.  If~$\beta \in \ZZ^d$,
the only localizations\/ $\NN A + \ZZ \mu$ capable of
containing~$\beta - m\alpha$ for all large (positive) integers~$m$ are
those for faces~$\mu$ containing~$\rho$.  In other words,
\begin{eqnarray*}
  \mu \in \nabla(\beta - m\alpha) \text{ for all } m \gg 0 &\implies&
  \mu \supseteq \rho.
\end{eqnarray*}
\end{lemma}
\begin{proof}
If~$\mu$ does not contain~$\rho$, then choose a linear functional that
is zero along~$\mu$ but positive on~$\alpha$.  This linear functional
remains negative on $\beta - m\alpha + \gamma$ for all $m \gg 0$ and
$\gamma \in \mu$, so that $\beta - m\alpha \not\in \NN A + \ZZ\mu$.%
\end{proof}

\begin{lemma} \label{l:rho}
Fix $\beta \in \ZZ^d$.  Suppose that $\rho$ is maximal among faces
of\/~$\NN A$ not in~$\nabla(\beta)$, but that $\rho$ is neither\/ $\NN
A$ nor a facet of\/~$\NN A$.  If~$\alpha \in \rho$ is a vector not
lying on any proper face of~$\rho$, then $\beta - m\alpha$ is an
exceptional degree for all large integers~$m$.%
\end{lemma}
\begin{proof}
Suppose that~$\mu$ contains~$\rho$.  Since $m\alpha \in \ZZ\mu$ for
all integers~$m$, we find that $\beta - m\alpha \in \NN A + \ZZ\mu$
for all integers~$m$ if and only if $\beta \in \NN A + \ZZ\mu$.
By~Lemma~\ref{l:m} we conclude that $\nabla(\beta - m\alpha)$ is, for
$m \gg 0$, the cocomplex of all faces strictly containing~$\rho$.  The
cohomology of such a cocomplex is the same as that of a sphere having
dimension $1+\dim(\rho)$: a~copy of~$\CC$ in dimension~$1+\dim(\rho)$
and zero elsewhere.  This cohomology is not in cohomological
degree~$d$ by the codimension hypothesis on~$\rho$.%
\end{proof}

\begin{lemma} \label{l:e}
Suppose that $\nabla$ is a cocomplex inside of a simplex of
dimension~$e$.  If the cohomology $H^j(\nabla)$ is nonzero for some $j
< e$, then there is a face~$\xi$ of codimension at least~$2$ inside
the simplex such that $\xi \not\in \nabla$ but $\mu \in \nabla$ for
all other faces~$\mu$ containing~$\xi$.%
\end{lemma}
\begin{proof}
The equivalent dual statement is easier to visualize: If $\Delta$ is a
simplicial complex inside of a simplex, and the reduced homology
$\tilde H_j(\Delta)$ is nonzero in some homological degree $j \geq 0$,
then there is a face~$\xi$ of dimension at least~$1$ in the simplex
such that $\xi \not\in \Delta$ but $\mu \in \Delta$ for every proper
face~$\mu$ of~$\xi$.  Equivalently, $\Delta$ has a minimal {\em
non\/}face of dimension at least~$1$.  This statement reduces easily
to the case where all vertices of the simplex lie in~$\Delta$, and in
that case one notes that {\em every}\/ nonface has dimension at
least~$1$, assuming $\Delta$ has at least two vertices.  But $\Delta$
has reduced homology in dimension $j \geq 0$, so it must have at least
two vertices and at least one nonface.%
\end{proof}

\begin{remark}
Lemma~\ref{l:e} fails immediately for polyhedral cocomplexes that are
not simplicial.  Two parallel edges of a square (plus the interior
cell of the square) form a polyhedral cocomplex that has cohomology of
dimension~$1$ in cohomological degree~$1$, but the only two maximal
nonfaces are the remaining two edges, of~codimension~$1$.
\end{remark}

\begin{thm}[\cite{mutsumi2}] \label{thm:exceptionalsetsimplex}
Suppose that the polytope\/ $\conv(A)$ is a simplex.  Then $\beta$ is
a rank-jumping parameter of~$A$ if and only if there exist faces
$\sigma$ and~$\tau$ of\/~$\NN A$, and an element\/ $\lambda \in
\CC\sigma\cap\CC\tau$, such that
\begin{eqnarray*}
  \lambda \in E_{\sigma}(\beta) \cap E_{\tau}(\beta) &\mathrm{but}&
  \lambda \not \in E_{\sigma \cap \tau}(\beta).
\end{eqnarray*}
\end{thm}

\begin{proof}[Proof of~Theorem~\ref{t:simplicial}]
We begin by showing that the exceptional set is contained in the slab
arrangement.  Let $\beta \in \CC^d$ be a rank-jumping parameter of
$A$, and pick $\sigma$, $\tau$ and $\lambda$ as in
Theorem~\ref{thm:exceptionalsetsimplex}.  For any $\alpha \in
\CC\sigma\cap\CC\tau$, the sum $\beta + \alpha$ is a rank-jumping
parameter, as can be seen by replacing~$\lambda$ with $\lambda +
\alpha$ in Theorem \ref{thm:exceptionalsetsimplex} and noting that
\begin{eqnarray} \label{alpha}
  \lambda+\alpha \in E_\tau(\beta+\alpha) &\iff& \lambda \in
  E_\tau(\beta).
\end{eqnarray}
Therefore we may (and do) assume that $\beta \in \ZZ^d$ and $\lambda =
0$.

Recall~\eqref{nabla}, which in particular implies that the set of
faces~$\mu$ satisfying $0 \in E_\mu(\beta)$ forms a cocomplex.  This
allows us to enlarge $\sigma$ and~$\tau$ so that $\sigma\cap\tau$ is
maximal among faces of~$\NN A$ outside~$\nabla(\beta)$, while still
satisfying Theorem~\ref{thm:exceptionalsetsimplex}.
Taking $\rho = \sigma\cap\tau$ in Lemma~\ref{l:rho},
we find that $\beta - m\alpha$ is an exceptional degree for all $m \gg
0$ and all choices of $\alpha$ interior to~$\rho$.  The slab
arrangement~$\ol E(A)$ therefore contains $\beta + \CC\rho$, and
hence~$\beta$.

Now suppose by Proposition~\ref{p:slab} that $\beta + \CC\rho$ is an
irreducible component of the slab arrangement $\ol E(A)$, where $\beta
\in \ZZ^d$ and $\rho$ is a face of~$\NN A$.  We wish to show that
$\beta + \CC\rho$ consists of rank-jumping parameters.  In fact, we
shall produce $\sigma$, $\tau$, and $\lambda$ as in
Theorem~\ref{thm:exceptionalsetsimplex} satisfying $\rho =
\sigma\cap\tau$, although we might harmlessly shift~$\beta$ by some
vector in~$\ZZ^d \cap \CC\rho$ first.

The component $\beta + \CC\rho$ is the closure of some slab parallel
to~$\rho$, which must (perhaps after shifting~$\beta$ by an element in
$\ZZ^d \cap \CC\rho$)
contain $\beta - m\alpha$ for an integer point~$\alpha$ interior
to~$\rho$ and all $m \gg 0$.  Replace $\beta$ by $\beta - m\alpha$ for
some fixed large choice of~$m$.  Lemma~\ref{l:m} implies that the
cocomplex $\nabla(\beta)$ is contained in the simplex consisting of
all faces of~$\NN A$ containing~$\rho$.  If~$\rho$ has dimension $d -
e - 1$, then this simplex satisfies the hypotheses of Lemma~\ref{l:e}.
Therefore we can find a face~$\xi$ containing~$\rho$ and of dimension
at most \mbox{$d-2$}, such that $\xi$ is a maximal nonface
of~$\nabla(\beta)$.  Applying Lemma~\ref{l:rho} to~$\xi$ instead
of~$\rho$, we find that the component $\beta + \CC\xi$ contains $\beta
+ \CC\rho$, and still lies inside~$\ol E(A)$.  {}From this we conclude
that $\rho = \xi$, because $\beta + \CC\rho$ is an irreducible
component of~$\ol E(A)$.

In summary, given that $\beta + \CC\rho$ is an irreducible component
of the slab arrangement~$\ol E(A)$, we have moved $\beta$ by an
element in~$\ZZ^d \cap \CC\rho$ so that
\begin{eqnarray*}
  \rho \not\in \nabla(\beta), &\text{but}& \mu \in \nabla(\beta)
  \text{ for all faces } \mu \text{ strictly containing } \rho.
\end{eqnarray*}
Moreover, we have shown that $\dim(\rho) \leq d-2$.  Therefore we can
pick two faces $\sigma$ and~$\tau$ strictly containing~$\rho$ and
satisfying $\rho = \sigma\cap\tau$.  For each $\lambda \in \CC\rho$,
we find that $\beta + \lambda$ is a rank-jumping parameter by
substituting $\lambda = 0$ and $\alpha = \lambda$ in \eqref{alpha},
then using~\eqref{nabla}, and finally applying
Theorem~\ref{thm:exceptionalsetsimplex}.%
\end{proof}

\subsection*{Acknowledgments}
The results in this note were obtained as part of a larger project
joint with Uli Walther, to whom we are very grateful.  We also thank
Bernd Sturmfels for his encouragement throughout this project.  The
Mathematical Sciences Research Institute (MSRI) in Berkeley,
California played a key role in the development of this project.  In
particular, the research reported in this article was completed while
both authors were postdoctoral fellows at MSRI, where LFM was
partially supported by a Postdoctoral Fellowship.  EM was supported by
a National Science Foundation Postdoctoral Research Fellowship.


\begin{thebibliography}{MMW04}

\bibitem[Ado94]{adolphson}
Alan Adolphson, \emph{Hypergeometric functions and rings generated by
  monomials}, Duke Math. J. \textbf{73} (1994), no.~2, 269--290.

\bibitem[BH93]{BH}
Winfried Bruns and J{\"u}rgen Herzog, \emph{Cohen-{M}acaulay rings},
  Cambridge University Press, Cambridge, 1993.

\bibitem[CDD99]{monomialcurves}
Eduardo Cattani, Carlos D'Andrea, and Alicia Dickenstein, \emph{The
  {${\mathcal{A}}$}-hyper\-geometric system associated with a monomial
  curve}, Duke Math. J. \textbf{99} (1999), no.~2, 179--207.

\bibitem[{CoC}]{cocoa}
{CoCoA}Team, \emph{{{\hbox{\rm C\kern-.13em o\kern-.07em C\kern-.13em
  o\kern-.15em A}}}: a system for doing computations in commutative
  algebra}, available at \textsf{http:/$\!$/cocoa.dima.unige.it/}.

\bibitem[GGZ87]{GGZ}
I.~M. Gelfand, M.~I. Graev, and A.~V. Zelevinsky, \emph{Holonomic
  systems of equations and series of hypergeometric type},
  Dokl. Akad. Nauk SSSR \textbf{295} (1987), no.~1, 14--19.

\bibitem[GPS01]{Singular}
G.-M. Greuel, G.~Pfister, and H.~Sch\"onemann, \emph{{\sc Singular}
  2.0}, \emph{A computer algebra system for polynomial computations},
  Centre for Computer Algebra, University of Kaiserslautern, 2001,
  \textsf{http:/$\!$/www.singular.uni-kl.de/}.

\bibitem[GS]{M2}
Daniel~R. Grayson and Michael~E. Stillman, \emph{Macaulay 2, a
  software system for research in algebraic geometry}, available at
  \textsf{http:/$\!$/www.math.uiuc.edu/Macaulay2/}.

\bibitem[GZK89]{GKZ}
I.~M. Gelfand, A.~V. Zelevinsky, and M.~M. Kapranov,
  \emph{Hypergeometric functions and toric varieties},
  Funktsional. Anal. i Prilozhen. \textbf{23} (1989), no.~2, 12--26.
  Erratum: Funktsional. Anal. i Prilozhen. \textbf{27} (1993), no.~4,
  91.


\bibitem[HM03]{injAlg}
David Helm and Ezra Miller, \emph{Algorithms for graded injective
  resolutions and local cohomology over semigroup rings}, J. Symbolic
  Computation, to appear, 2004.  \textsf{arXiv:math.AG.0309256}

\bibitem[MMW04]{ranks}
Laura~Felicia Matusevich, Ezra Miller, and Uli Walther,
  \emph{Homological methods for hypergeometric families}, in
  preparation, 2004.

\bibitem[Sai02]{mutsumi2}
Mutsumi Saito, \emph{{Logarithm-free {$A$}-hypergeometric series}},
  Duke Math. J. \textbf{115} (2002), no.~1, 53--73.

\bibitem[ST98]{StuTaka}
Bernd Sturmfels and Nobuki Takayama, \emph{Gr\"obner bases and
  hypergeometric functions}, Gr\"obner bases and applications (Linz,
  1998), London Math. Soc. Lecture Note Ser., vol. 251, Cambridge
  Univ. Press, Cambridge, 1998, pp.~246--258.

\end{thebibliography}


\end{document}